\documentclass[10pt,thmsa]{article}%
\usepackage{amsmath}
\usepackage{amsfonts}
\usepackage{amssymb}
\usepackage{graphicx}%
\newtheorem{theorem}{Theorem}[section]

\newtheorem{lemma}[theorem]{Lemma}

\newtheorem{proposition}[theorem]{Proposition}

\newenvironment{proof}[1][Proof]{\noindent\textbf{#1.} }{\ \rule{0.5em}{0.5em}}

\begin{document}
\title{Rigidity of negatively curved geodesic orbit Finsler spaces
\thanks{Supported by NSFC (no. 11271216, 11271198, 51535008), SRFDP of China, Science and Technology Development Fund for Universities and Colleges in Tianjin (no. 20141005), and Doctor fund of Tianjin Normal University (no. 52XB1305)}}
\author{Ming Xu$^1$ and Shaoqiang Deng$^2$\thanks{S. Deng is the corresponding author} \\
\\
$^1$College of Mathematics,
Tianjin Normal University\\
 Tianjin 300387, P. R. China\\
 Email: mgmgmgxu@163.com.\\
 \\
$^2$School of Mathematical Sciences and LPMC,
Nankai University\\
Tianjin 300071, P. R. China\\
E-mail: dengsq@nankai.edu.cn}
\date{}
\maketitle

\noindent\textbf{Abstract}

We prove some rigidity results on geodesic orbit Finsler spaces with non-positive  curvature. In particular, we show that a geodesic Finsler space with strictly negative flag curavture must be a non-compact Riemannian symmetric space of rank one.

\noindent\textbf{Key Words}: Finsler space; geodesic orbit space; flag curvature.

%

\section{Introduction}
The main goal of this short note is to prove some rigidity results on  geodesic orbit Finsler spaces (simply g.o. Finsler spaces) with negative curvature.
Recall that a homogeneous Finsler space $(G/H,F)$ is called a
{\it $G$-geodesic orbit Finsler space}, if any geodesic of $(G/H,F)$ is {\it $G$-homogeneous}, that is,  it is the orbit of a  one-parameter subgroup
$g_t$ in $G$ generated by a vector $v\in\mathrm{Lie}(G)=\mathfrak{g}$. This definition is a natural generalization of Riemannian geodesic orbit space, which has been
extensively studied, see for example \cite{AA, GD}. Finslerian g.o. spaces and homogeneous normal Finsler spaces have also recently been studied in \cite{XD2014} and \cite{YD}. The research shows that
 there exists many
examples of compact Riemannian or non-Riemannian g.o. spaces, e.g.,  Clifford-Wolf  homogeneous spaces, normal homogeneous spaces, generalized normal homogeneous spaces (or $\delta$-homogeneous spaces), etc. Some noncompact examples have also constructed.

Our rigidity results are relevant to curvatures in Finsler geometry, where flag curvature is one of the most important geometric quantities, and most of the research in this field is focused on issues related to it. For example,  the classification of Finsler spaces with certain flag curvature conditions is a very interesting problem. Recently, we have made some progress in this consideration, see \cite{DX2014, XD2014, XDHH2014}.
 For example,  a complete classification of
normal homogeneous Finsler spaces with positive flag curvature  is established, which generalizes the results of
M. Berger \cite{Ber61}. In this process we found that many positively curved normal homogeneous Riemannian manifolds can be endowed with invariant positively curved non-Riemannian normal Finsler metrics.  One may expect that may also exist non-Riemannian negatively curved examples. However, Contrary to the above fact, in the negative case we get the following rigidity theorem.

\begin{theorem} \label{mainthm}
Any connected negatively curved g.o. Finsler space must be a compact
Riemannian symmetric space of rank one.
\end{theorem}

 This fact is interesting in that any homogeneous normal Finsler space must be a g.o. space \cite{XD2014}. We have also made some efforts to generalize the result to   g.o. Finsler spaces  with negative  Ricci curvature. However,
 the problem becomes  rather hard, and it is still open. Nevertheless,
the same argument can also be applied to a slightly general case to get the following  rigidity theorem.

\begin{theorem}\label{mainthm-2}
Any g.o. Finsler space with negative Ricci scalar and non-positive flag curvature must be a symmetric space of non-compact type.
\end{theorem}

Note that a symmetric Finsler space   must be Berwald, but  can be non-Riemannian if its rank is larger than one \cite{Deng}.

Another goal of this paper is to  explore the general  properties of
negatively curved homogeneous Finsler metrics.  We will show that, any connected negatively curved homogeneous Finsler space can be viewed as a solvemanifold endowed with a left invariant Finsler metric. This fact motivates us to consider the classification of left invariant  Finsler metrics on a solvable Lie group  with negative curvature. Applying the homogeneous flag curvature formula of L. Huang \cite{Huang-1}, we   prove

\begin{theorem}\label{mainthm-3}
Let $G$ be a  simply connected solvable Lie group such that   $[\mathfrak{g},\mathfrak{g}]$ is Abelian and nonzero, where $\mathfrak g$ is the Lie algebra of $G$. Then it admits a negatively curved left invariant Finsler metric if and only if the following conditions are satisfied:
\begin{description}
\item{\rm (1)} $\dim\mathfrak{g}=\dim[\mathfrak{g},\mathfrak{g}]+1$;
\item{\rm (2)} There exists a vector $u\in\mathfrak{g}$, such that all the eigenvalues of
$\mathrm{ad}(u)|_{[\mathfrak{g},\mathfrak{g}]}: [\mathfrak{g},\mathfrak{g}]
\rightarrow[\mathfrak{g},\mathfrak{g}]$,
have positive real parts.
\end{description}
\end{theorem}

Inspired by \cite{He74}, we conjecture that Theorem \ref{mainthm-3} is  valid for
 general negatively curved homogeneous Finsler spaces.


\medskip
\noindent {\bf Acknowledgement.} We are grateful to  J. A. Wolf and L. Huang for helpful discussions.
\section{A  homogeneous flag curvature formula}

On a Finsler space $(M,F)$, the Riemann curvature
$R_y^F=R_k^i(y)\partial_{x^i}\otimes dx^k: T_x M\rightarrow T_xM$ can be similarly defined as in
Riemann geometry, either by the structure equation of the Chern connection,
or the Jacobi field equation of the variation of geodesics \cite{Shen}, and it
 can be used to define  flag curvature and  Ricci curvature. Let $y\in T_xM$
be a nonzero tangent vector and $\mathbf{P}$ a tangent plane in $T_xM$
containing $y$, and suppose $\mathbf{P}$ is linearly spanned by $y$ and $v$. Then the flag
curvature of the triple $(x,y,y
\wedge v)$ or $(x,y,\mathbf{P})$ is defined as
\begin{equation*}
K^F(x,y,y\wedge v)=K^F(x,y,\mathbf{P})=
\frac{\langle R_y v,v\rangle_y}{\langle y,y\rangle_y\langle v,v\rangle_y
-\langle y,v\rangle_y^2},
\end{equation*}
and the Ricci scalar of the nonzero vector $y\in T_xM$ is
$$\mathrm{Ric}^F(x,y)=
F(x,y)^2\sum_{i=1}^{n-1}K^F(x,y,y\wedge e_i)$$
where $\{e_1,\ldots,e_{n-1},y/F(y)\}$ is a $g_y^F$-orthonormal basis of $T_xM$.


The following flag curvature formula, which is very convenient in homogeneous Finsler geometry, can be found in \cite{XDHH2014}.

\begin{theorem} \label{flag-curvature-formula-thm}
Let $(G/H,F)$ be a connected homogeneous Finsler space, and $\mathfrak{g}=\mathfrak{h}+\mathfrak{m}$ be an $\mathrm{Ad}(H)$-invariant
decomposition for $G/H$. Then for any linearly independent commutative pair
$u$ and $v$ in $\mathfrak{m}$ satisfying
$
\langle[u,\mathfrak{m}],u\rangle^F_u=0
$,
we have
\begin{equation*}
K^F(o,u,u\wedge v)=\frac{\langle U(u,v),U(u,v)\rangle_u^F}
{\langle u,u\rangle_u^F \langle v,v\rangle_u^F-
{\langle u,v\rangle_u^F}\langle u,v\rangle_u^F},
\end{equation*}
where $U$ is a bilinear map from $\mathfrak{m}\times \mathfrak{m}$ to $\mathfrak{m}$ defined by
\begin{equation*}
\langle U(u,v),w\rangle_u^F=\frac{1}{2}(\langle[w,u]_\mathfrak{m},v\rangle_u^F
+\langle[w,v]_\mathfrak{m},u\rangle_u^F), \,\,w\in\mathfrak{m},
\end{equation*}
here $[\cdot,\cdot]_\mathfrak{m}=\mathrm{pr}_\mathfrak{m}\circ[\cdot,\cdot]$ and $\mathrm{pr}_\mathfrak{m}$ is the projection
with respect to the given $\mathrm{Ad}(H)$-invariant decomposition.
\end{theorem}

\section{Proofs of Theorem \ref{mainthm} and Theorem \ref{mainthm-2}}

\noindent\textbf{Proof of Theorem \ref{mainthm}}.\quad
Let $(M,F)$ be a connected negatively curved geodesic orbit Finsler space.
Then we can write $M$ as $M=G/H$, where $G=I_0(M,F)$ and $H$ is the compact isotropy subgroup at $o=eH$.

 We first prove that $G$ is a semisimple Lie group. Assume conversely that it is not.
Then we can find an Abelian ideal $\mathfrak{a}$ of $\mathfrak{g}$ such that $\mathfrak{a}\cap\mathfrak{h}=\{0\}$, where $\mathfrak{h}=\mathrm{Lie}\,H$. Hence we can find an
$\mathrm{Ad}(H)$-invariant decomposition
$\mathfrak{g}=\mathfrak{h}+\mathfrak{m}$ such that $\mathfrak{a}\subset\mathfrak{m}$.

\begin{lemma}\label{lemma-1}Keep the above notation and assumptions.
Then any vector in
$\mathfrak{a}$ defines
a Killing vector field of constant length of $(G/H,F)$.
\end{lemma}

\begin{proof}
We first prove that $o=eH$ is the critical point of the function
$F^2(V(\cdot))$, where $V$ is the Killing vector field defined by a fixed
 nonzero vector $v\in\mathfrak{a}$. This is equivalent to
$\langle v,[v,\mathfrak{g}]_\mathfrak{m}\rangle_v^F=0$. Since
$\langle v,[v,\mathfrak{h}]\rangle_v^F=0$, we only need to show that
$\langle v,[v,\mathfrak{m}]_\mathfrak{m}\rangle_v^F=0$. By Lemma 3.3 in \cite{YD},
there exists a vector $v'\in\mathfrak{h}$, such that
$\langle v,[v,u]_\mathfrak{m}\rangle_v^F=
\langle v,[v',u]\rangle_v^F, \forall u\in \mathfrak{m}.$
On the other hand, by Theorem 3.1 of \cite{DH04},
we have
$$\langle v,[v',u]\rangle_v^F=\langle[v,v'],u\rangle_v^F+2C^F_v(v,u,[v,v'])
=\langle[v,v'],u\rangle_v^F.$$
If $u$ is contained in the $g_v^F$-orthogonal complement of $\mathfrak{a}$
in $\mathfrak{m}$, then $\langle[v,v'],u\rangle_v^F
\subset\langle \mathfrak{a},v'\rangle_v^F=0.$  In case $u\in \mathfrak{a}$, we also have $\langle v,[v,u]_\mathfrak{m}\rangle_v^F=0$.
Therefore  we always have $\langle v,[v,\mathfrak{m}]_\mathfrak{m}\rangle_v^F=0$.
Thus $o=eH$ is a critical point of the function $F^2(V(\cdot))$.

Now given $g\in G$, $\mathrm{Ad}(g^{-1})v$ is a nonzero
vector in $\mathfrak{a}$, which defines another Killing vector field $V'$
of $(G/H,F)$. Moreover, we have $F(V'(\cdot))=F(V(g\cdot))$, and $F(V'(e))=F(V(g))\neq 0$. Since $o=eH$
is a critical point of $F^2(V'(\cdot))$,
$x=gH$ is also a critical point of $F^2(V(\cdot))$. Consequently each point is a critical point of $F^2(V(\cdot))$, that is,  $V$ is a Killing vector field of constant
length.
\end{proof}

By the above lemma, both the flag curvature $K^F(o,y,\mathbf{P})$ and the Ricci scalar
$\mathrm{Ric}^F(o,y)$ are non-negative for nonzero $y\in\mathfrak{a}$
(see the proof of  Theorem 5.1 in \cite{XD2014} for more details).
This is a contradiction. Thus $G$ is semisimple.

Next we prove that $H$ is a maximal compact subgroup of $G$.
Assume conversely that it is not. Then $H$ is contained in a connected
maximal subgroup $K$ with $\dim H<\dim K$. Correspondingly we have the
Iwasawa decomposition $G=NAK$, where the solvable subgroup $NA$ is diffeomorphic to a Euclidean space. Then $G/H$ is diffeomorphic to
$NA\times K/H$, hence it is homotopic to $K/H$. By the main theorem in \cite{DH07}, $G/H$ is diffeomorphic to a Euclidean space, but the compact
homogeneous space $K/H$ can not be topologically trivial. This  contradiction proves our assertion.

In summarizing, $H$ is a maximal compact subgroup of the semisimple Lie group $G$, hence
$G/H$ is a simply connected symmetric homogeneous space of non-compact type.
The homogeneous metric $F$ on $G/H$ is Berwald, so it is negatively curved if and only if the
Riemannian symmetric homogeneous space $G/H$ is negatively curved, i.e., $G/H$ is of rank 1 \cite{Deng}. However, in an irreducible Riemannian symmetric spaces of rank 1, the isotropy group $H$ acts transitively on the unit sphere of the tangent space at the origin. So $F$ must be a positive multiple of the canonical  Riemannian metric.
This completes the proof of the theorem.

\noindent\textbf{Proof of Theorem
\ref{mainthm-2}}.\quad Note that the arguments in the proof of Theorem
\ref{mainthm} are still valid under the assumption that $(G/H,F)$ has non-positive flag curvature and negative Ricci scalar.
So we can similarly show that $G/H$ is a simply connected symmetric space, but not necessarily of rank $1$. Then the homogeneous Finsler metric $F$ must be Berwald, and $F$ can be  non-Riemannian (see \cite{Deng}).

\section{Negatively curved solvemanifolds}

We first show that any connected negatively curved homogeneous Finsler space can be viewed as a solvemanifold endowed with a left invariant Finsler metric.
Let $G/H$ be a connected negatively curved homogeneous Finsler space, where
$G$ is connected and $H$ is the compact isotropy subgroup at $o=eH$.
By \cite{DH07}, $G/H$
is simply connected. Using Theorem 3.13 in \cite{Iw49},
we can also  prove that $H$ is a maximal compact subgroup of $G$ using similar
argument as in the proof of Theorem \ref{mainthm}, and this proves our assertion.
We remark here that, recently, Wang and Li \cite{WL} make some arguments to show that any connected negatively curved homogeneous Finsler space $(M,F)$ admits a closed solvable transitive subgroup $G$ of $I_0(M,F)$, based on the results of Wolf in \cite{Wo1964}.

{\bf Proof of Theorem \ref{mainthm-3}}.
If $G$ satisfies (1) and (2) in the theorem, then by \cite{He74}
it admits a left invariant negatively curved Riemannian metric.

Assume that $G$ is a simply connected solvable Lie group, with
$\mathrm{Lie}(G)=\mathfrak{g}$ satisfying $[\mathfrak{g},\mathfrak{g}]=0$, and $F$  is a left invariant
negatively curved Finsler metric on $G$.
Then (1) follows from a unpublished proposition of L. Huang. With his
permission, we quote it here.

\begin{proposition}
Let $F$ be a left invariant negatively curved Finsler metric on
a simply connected solvable Lie group $G$. Then
$\mathfrak{g}=\mathrm{Lie}(G)$ must satisfy
$\dim\mathfrak{g}=\dim[\mathfrak{g},\mathfrak{g}]$+1.
\end{proposition}

\begin{proof}
 Since $G$ (and $\mathfrak{g}$) is solvable, we have $\dim[\mathfrak{g},\mathfrak{g}]<\dim\mathfrak{g}$.
Thus there exists a nonzero vector $u\in\mathfrak{g}\backslash[\mathfrak{g},\mathfrak{g}]$,
such that $\langle u,[\mathfrak{g},\mathfrak{g}]\rangle_u^F=0$.
Then a direct calculation shows that $\eta(u)=0$.
If there exists a nonzero vector
$v\in\mathfrak{g}$ such that $u$ and $v$ are  commutative and linearly independent, then we have $\langle u,[u,\mathfrak{g}]\rangle_u^F=0$. Hence   by
Theorem \ref{flag-curvature-formula-thm},
$K^F(o,u,u\wedge v)\geq 0$.
Thus $\mathrm{Ad}(u):
\mathfrak{g}\rightarrow[\mathfrak{g},\mathfrak{g}]$
has only a one-dimensional kernel $\mathbb{R}u$, i.e.,
$\dim\mathfrak{g}=\dim[\mathfrak{g},\mathfrak{g}]+1$.
\end{proof}

Now we prove (2). If $\dim[\mathfrak{g},\mathfrak{g}]=1$, then either (2) is
satisfied or $\mathfrak{g}$ is Abelian. But in the second case, $(G,F)$
is flat.

Assume $\dim[\mathfrak{g},\mathfrak{g}]>1$.  We assert that there does not exist a pair of nonzero vectors
$u\in[\mathfrak{g},\mathfrak{g}]$ and
$u'\in\mathfrak{g}\backslash[\mathfrak{g},\mathfrak{g}]$ such that
$\langle u,[u,u']\rangle_u^F=0$. Otherwise, we have $\langle u,[u,\mathfrak{g}]\rangle_u^F=0$, since
$[\mathfrak{g},\mathfrak{g}]=0$. We can then find a vector $v\in[\mathfrak{g},\mathfrak{g}]$ such that
$u$ and $v$ are linearly independent. Then by Theorem \ref{flag-curvature-formula-thm}, $K^F(o,u,u\wedge v)\geq 0$.
This is a contradiction.

Consequently we can find a vector $u'\in\mathfrak{g}\backslash[\mathfrak{g},\mathfrak{g}]$,
such that $\langle u,[u',u]\rangle_u^F>0$
for any nonzero $u\in[\mathfrak{g},\mathfrak{g}]$.
Let $g_t$ be the one-parameter subgroup of $G$ generated by $u'$.  Then
\begin{equation}
\frac{d}{dt}F(\mathrm{Ad}(g_t)u)|_{t=0}>0\label{0},
\end{equation}
 for any nonzero $u\in[\mathfrak{g},\mathfrak{g}]$.
 Now suppose that the restriction of $\mathrm{Ad}(u')$ to $[\mathfrak{g},\mathfrak{g}]$ has an eigenvalue whose real part is non-positive. If this eigenvalue is real,
then $\langle u,[u',u]\rangle_u^F\leq 0$ for the nonzero eigenvalue $u$.
On the other hand, if this eigenvalue is not real, then we can find a linearly independent pair $u_1$ and
$u_2$, such that all $\mathrm{Ad}(g_t)$ with $t>0$ rotate two dimensional subspace $\mathbf{V}=\mathbb{R}u_1+\mathbb{R}u_2$, or shrink it at the same time as rotate $\mathbf{V}$. This is a contradiction with (\ref{0}) restricted
to $\mathbf{V}$. Thus  all  the eigenvalues of $\mathrm{Ad}(u')|_{[\mathfrak{g},\mathfrak{g}]}$ have positive real parts, proving (2).

\end{document}